\documentstyle[12pt]{article}
 \newcommand{\beq}{\begin{equation}}
\newcommand{\eeq}{\end{equation}}
 \newcommand{\beth}{\begin{theo}}
\newcommand{\eth}{\end{theo}}

\newcommand{\Subset}{\subset\subset}

\newcommand{\Psux}{{ \Psi_{u,x}}}

\newcommand{\psux}{{ \psi_{u,x}}}
\newcommand{\nux}{\nu(u,x)}
\newcommand{\nui}{\nu(u,\infty)}
\newcommand{\nuai}{\nu(u,a,\infty)}

\newcommand{\tmux}{{\cal T}_{m,x} u}
\newcommand{\tsmux}{{\cal T}_{m_s,x} u}
\newcommand{\cL}{{Log}}
\newcommand{\cE}{{Exp}}
\newcommand{\okn}{1\le k\le n}
\newcommand{\LL}{{\cal L}}
\newcommand{\mucn}{M(u;\Cn)}
\newcommand{\TT}{T}
\newcommand{\II}{{\cal I}}
\newcommand{\Ed}{{\bf 1}}

\newcommand{\Znp}{{ \bf Z}_+^n}
\newcommand{\Rn}{{ \bf R}^n}

\newcommand{\Rnp}{{ \bf R}_+^n}
\newcommand{\Cn}{{{ \bf C}^n}}

\newcommand{\CN}{{ \bf C}^N}
\newcommand{\Cstar}{{{ \bf C}^{*n}}}

\newtheorem{theo}{Theorem}

\newtheorem{prop}{Proposition}

\begin{document}

\begin{center}
 {\bf \large  INDICATORS FOR PLURISUBHARMONIC
FUNCTIONS}
\end{center}
\begin{center}{\bf \large OF LOGARITHMIC GROWTH} \end{center}

\vskip0.5cm
\begin{center}
{\bf Alexander RASHKOVSKII  }
\end{center}

\vskip1cm

  {\small {\sc Abstract.}  A notion of indicator for a
plurisubharmonic function $u$ of logarithmic growth in $\Cn$ is
introduced and studied. It is applied to evaluation of the total
Monge-Amp\`ere measure  $(dd^cu)^n(\Cn)$. Upper bounds for 
the measure are obtained in terms of growth characteristics of $u$. 
 When $u=\log|f|$ for a polynomial mapping $f$
with isolated zeros, the indicator generates the Newton polyhedron of $f$
whose volume bounds the number of the zeros. 
}

\vskip 0.5cm

\section{Introduction}

We consider plurisubharmonic functions $u$
of logarithmic growth in $\Cn$, i.e. satisfying the relation
\beq
u(z)\le C_1\log^+|z|+C_2
\label{eq:1}
\eeq
with some constants $C_j=C_j(u)\ge 0$. The class of such functions
will be denoted by $\LL(\Cn)$ or simply by $\LL$. (It is worth
mentioning that in the literature the notation $\LL$ is used
sometimes for the class of functions satisfying (\ref{eq:1}) with
$C_1=1$; for our purposes we need to consider the whole class
of functions of logarithmic growth, and denoting it by $\LL$ we
follow, for example, \cite{keyLe4}, \cite{keyLe5}.) It is an
important class containing, in particular, functions of the form
$\log|P|$ with polynomial mappings $P:\Cn\to\CN$. Various results
concerning the functions of logarithmic growth can be found in
\cite{keyLe1}-\cite{keyLe5}, \cite{keySi1}, \cite{keyKl}, see also
the references in \cite{keyKl} and \cite{keyKis5}. For general
properties of plurisubharmonic functions and the complex
Monge-Amp\`ere operators, we refer the
reader to \cite{keyLe2}, \cite{keyLeG}, \cite{keyKl}, and
\cite{keyD1}.

A remarkable property of functions $u\in\LL$ is finiteness of their
total Monge-Amp\`ere measures
$$
M(u;\Cn)=\int_\Cn(dd^cu)^n<\infty
$$
as long as $(dd^cu)^n$ is well defined on the whole $\Cn$; we use
the notation $d=\partial + \bar\partial,\ d^c= ( \partial
-\bar\partial)/2\pi i$. Moreover, the total mass is tied strongly to
the growth of the function. For example, if
\beq
\log^+|z|+c\le u(z)\le \log^+|z|+C,
\label{eq:2}
\eeq
then $\mucn=M(\log|z|;\Cn)=1$. The objectives for the present paper
is to study $\mucn$ when no regularity condition on $u$ like
(\ref{eq:2}) is assumed. In case of $u=\log|P|$ with $P:\Cn\to\CN$
a polynomial mapping with isolated zeros, $\mucn$ equals (if $N=n$) or
dominates (if $N>n$) the number of the zeros counted with their
multiplicities.

If $u=v$ near the boundary of a bounded pseudoconvex domain $\Omega$,
then
$$
\int_\Omega (dd^cu)^n=\int_\Omega (dd^cv)^n,
$$
so the total measure $\mucn$ is determined by the asymptotic
behavior of $u$ at infinity. For its evaluation we thus need precise
characteristics of the behavior. The basic one is the {\it
logarithmic type}
\beq
\sigma(u)=\limsup_{z\to\infty} {u(z)\over\log|z|}.
\label{eq:type}
\eeq
Another known characteristic is the {\it logarithmic multitype}
$(\sigma_1(u),\ldots,\sigma_n(u))$ \cite{keyLe5}:
\beq
\sigma_1(u)=\sup\,\{\tilde\sigma_1(u;z'): z'\in{\bf C}^{n-1}\}
\label{eq:mtype}
\eeq
where $\tilde\sigma_1(u;z')$ is the logarithmic type of the function
$u_{1,z'}(z_1)=u(z_1,z')\in\LL({\bf C})$ with $z'\in{\bf C}^{n-1}$
fixed, and similarly for $\sigma_2(u),\ldots,\sigma_n(u)$.
For example, if $P$ is a polynomial of degree $d_k$ in $z_k$, then
$\sigma_k(\log|P|)=d_k$.

Due to the certain symmetry between the behavior of $u\in\LL$ at
infinity and the local behavior of a plurisubharmonic function at a
fixed point of its logarithmic singularity, the type $\sigma(u)$ can
be regarded as the
Lelong number of $u$ at infinity:
\beq
\sigma(u)=\nui.
\label{eq:lel}
\eeq
One can also consider the directional Lelong numbers at infinity with
respect to directions $a=(a_1,\ldots,a_n)\in\Rnp$:
\beq
\nuai=\limsup_{z\to\infty} {u(z)\over S_a(z)},
\label{eq:dir}
\eeq
where
\beq
S_a(z)=\sup_k\,a_k^{-1}\log|z_k|.
\label{eq:Sa}
\eeq

In \cite{keyLeR}, the residual Monge-Amp\`ere measure of a
plurisubharmonic function $u$ at a point $x\in\Cn$,
$(dd^cu)^n|_{\{x\}}$, was studied by means of the local indicator of
$u$ at $x$. Using the same approach, we introduce here a notion of
the indicator of $u\in\LL$:
$$
\Psux(y)=\lim_{R\to +\infty}R^{-1}\sup\{u(z):\: |z_k-x_k|\le
               |y_k|^R,\ \okn\}.
$$
It is a plurisubharmonic function of the class $\LL$ which is the
(unique) logarithmic tangent to $u$ at $x$, i.e. the weak limit in
$L_{loc}^1(\Cn)$ of the functions $m^{-1} u(x_1+y_1^m,\ldots,
x_n+y_n^m)$ as $m\to\infty$ (Theorem \ref{theo:tangent}).
The above characteristics of $u$ can be easily expressed in terms of
its indicator (see Proposition \ref{prop:char}); 
moreover, it controls the behavior of $u$ in the whole $\Cn$
(Theorem \ref{theo:maj}):
$$
u(z)\le\Psux(x-x)+C_x\quad\forall z\in\Cn.
$$

If $(dd^cu)^n$ is defined on $\Cn$, the indicator also controls the
total Monge-Amp\`ere mass of $u$ (Theorem \ref{theo:comp2}):
\beq
\mucn\le M(\Psux;\Cn).
\label{eq:Ibound}
\eeq
Since $\Psux(y)=\Psux(|y_1|,\ldots,|y_n|)$, the evaluation of its
mass is much more easy than that for the original function $u$. It gives
us, in particular, the bounds
$$
\mucn\le{[\nuai]^n\over a_1\ldots a_n}\quad\forall a\in\Rnp
$$
(Theorem \ref{theo:ocenka}) and
$$
\mucn\le n!\,\sigma_1(u)\ldots\sigma_n(u)
$$
(Theorem \ref{theo:Dyson}). A particular case of the latter
result (when $u$ is the logarithm of modulus of an
equidimensional polynomial mapping with isolated zeros
of regular multiplicities) was obtained in \cite{keyLe6}.

In Theorems \ref{theo:volume} and \ref{theo:Kou} we give a
geometric description for the mass of an indicator. Denote
$\psux(t)=\Psux(e^{t_1},\ldots,e^{t_n})$, $t=(t_1,\ldots,t_n)\in\Rn$,
and
$$
\Theta_{u,x}=\{a\in\Rn:\: \langle a,t\rangle \le\psux^+(t)\ \forall
t\in\Rn\}.
$$
Then
\beq
M(\Psux;\Cn)=n!\,Vol(\Theta_{u,x}).
\label{eq:IKou}
\eeq
When $u=\log|P|$ with $P$ a polynomial mapping, the set $\Theta_{u,0}$ is
the Newton polyhedron for $P$ at infinity (see, for example,
\cite{keyKou1}), i.e. the convex hull of the set
$\omega_0\cup\{0\}$,
$$ \omega_0=\{s\in\Znp:\: \sum_j \left|
\frac{\partial^s P_j}{\partial z^s} (0)\right|\neq 0\},
$$
and so the right-hand side of (\ref{eq:IKou}) is the Newton number of
$P$ at infinity. Therefore, a hard result due to A.G.~Kouchnirenko
on the number of zeros of an equidimensional polynomial mapping
\cite{keyKou} follows directly from (\ref{eq:Ibound}) and
(\ref{eq:IKou}).


\section{Indicators as growth characteristics}

Let $u$ be a plurisubharmonic function in $\Cn$. Given $x\in\Cn$
and $t\in\Rn$, denote by $g(u,x,t)$ the mean value of $u$ over the set
$\TT_t(x)=\{z\in\Cn:\;|z_k-x_k|=e^{t_k},\ \okn\}$, and by $g'(u,x,t)$
the maximum of $u$ on $\TT_t(x)$.

\begin{prop}
\label{prop:psi}
Let $u\in\LL,\ x\in\Cn$. Then for every $t\in\Rn$ the following
limits exist and coincide:
$$
\lim_{R\to +\infty} R^{-1} g(u,x,Rt)=
\lim_{R\to +\infty} R^{-1} g'(u,x,Rt)=: \psux(t)<\infty.
$$
Moreover, if $g'(u,x,0)\le 0$, the common limit $\psux(t)$ is obtained by
the increasing values.
\end{prop}

{\it Proof.} For $x\in\Cn$ and $t\in\Rn$ fixed, the function
$f(R):=g(u,x,Rt)$ is convex on ${\bf R}$ and has the bound
$f(R)\le C_1 R+C_2\ \forall R>0$ with some $C_1,C_2>0$. Therefore,
for all $R_0\in{\bf R}$, the ratio
\beq
f(R)-f(R_0)\over R-R_0
\label{eq:ratio}
\eeq
is increasing in $R>R_0$ and bounded and thus has a limit as
$R\to+\infty$. It implies the existence of
$\hat g(u,x,t)=\lim_{R\to +\infty} R^{-1} g(u,x,Rt)$.
In the same way we get the value
$\hat g'(u,x,t)=\lim_{R\to +\infty} R^{-1} g'(u,x,Rt)$.
Evidently, $\hat g(u,x,t)\le \hat g'(u,x,t)$, and the standard
arguments using Harnack's inequality give us
$\hat g(u,x,t)= \hat g'(u,x,t)$.
The last statement of the proposition follows from the increasing of
(\ref{eq:ratio}) with $R_0=0$.

\medskip
Now we proceed, as in \cite{keyLeR}, to a plurisubharmonic characteristic
of growth for $u\in\LL$. Denote
$\Cstar=\{z\in\Cn: z_1\ldots z_n\neq 0\}$.
The mappings $\cL:\Cstar\to\Rn$ and $\cE:\Rn\to \Cstar$
are defined as $\cL(z)=(\log|z_1|,\ldots,
\log|z_n|)$ and $\cE(t)=(\exp t_1,\ldots, \exp t_n)$, respectively.
Let $\LL^c$ be the subclass of $\LL$ formed by $n$-circled
plurisubharmonic functions $u$, i.e. $u(z)=u(|z_1|,\ldots,|z_n|)$.
By $\LL(\Rn)$ we denote the class of functions $\varphi(t)$, $t\in\Rn$,
which are convex in $t$, increasing in each $t_k$ and such that there
exists a limit $\lim_{T\to +\infty}T^{-1}\varphi(T,\ldots,T)<\infty$.

The mappings $\cE$ and $\cL$ generate an isomorphism between the cones
$\LL^c$ and $\LL(\Rn)$ (\cite{keyLe5}, Th.~1):
$u\in \LL^c \iff \cE^*u\in
\LL(\Rn)$, $h\in \LL(\Rn) \iff \cL^*h$ extends to a (unique) function from
the class $\LL^c$. Given $u\in\LL^c$, the function $\cE^*u$ will be
referred to as the {\it convex image} of $u$.

If $h=\cE^*u\in\LL(\Rn)$ satisfies the homogeneity condition
$$
h(ct)=ch(t)\quad\forall c>0,\ \forall t\in\Rn,
$$
the function $u$ will be called {\it an indicator}. We denote the
collection of all indicators by $\II$. It is easy to see that any
indicator $\Psi$ satisfies $\Psi\le 0$ in the unit polydisk
$$
D=\{z\in\Cn:\; |z_k|<1,\ \okn\}
$$
and $\Psi>0$ in
$$
D^{-1}=\{z\in\Cn:\; |z_k|>1,\ \okn\}
$$
if $\Psi\not\equiv 0$.

Clearly, the function $\psux$ defined in Proposition \ref{prop:psi}
belongs to the class $\LL(\Rn)$, so $\cL^*\psux$ extends to a function
$\Psux\in\LL^c$:
$$
\Psux(y)=\psux(\log|y_1|,\ldots,\log|y_n|),\ y\in\Cstar.
$$
Moreover, $\Psux\in\II$. We will call it {\it the indicator of
$u\in\LL$ at $x$}.

The restriction of $\Psux$ to the polydisk $D$ coincides with the local
indicator of $u$ at $x$ introduced in \cite{keyLeR}. In particular,
$\Psux\equiv 0$ in $D$ if and only if the Lelong number of $u$ at $x$
equals $0$. Besides, the directional Lelong numbers of $\Psux$ at $0$
are the same as those of $u$ at $x$.

\begin{prop}
\label{prop:indind}
Let $\Phi\in\II$, then
\begin{enumerate}
\item[(a)]
$\Phi$ is continuous as a function $\Cn\to{\bf R}\cup\{-\infty\}$;
\item[(b)]
$\Psi_{\Phi,x}(y)=\Phi(\tilde y)$ where $\tilde y_k=\sup\,\{|y_k|,1\}$
if $x_k\neq 0$, and $\tilde y_k=y_k$ otherwise.
\end{enumerate}
\end{prop}

{\it Proof.}
{\it (a)} Since $\cE^*\Phi\in C(\Rn)$, $\Phi\in C(\Cstar)$. Its continuity
on $\Cn$ can be shown by induction in $n$. Let it be already proved for
$n\le l$ (the case $n=1$ is obvious). Consider any point $z^0\in{\bf
C}^{l+1}$ with $z_j^0=0$ for some $j$. If $\Phi(z^0)=-\infty$, then
$\Phi(z^s)\to -\infty$ for every sequence $z^s\to z^0$. If
$\Phi(z^0)>-\infty$, consider the projections $\tilde z^s$ of
$z^s\to z^0$ to the subspace $L_j=\{z\in{\bf C}^{l+1}: z_j=0\}$:
$\tilde z_j^s=0$ and $\tilde z_m^s=z_m^s\ \forall m\neq j$. Since
$\Phi|_{L_j}\not\equiv -\infty$, the induction assumption implies
$\Phi(\tilde z^s)\to\Phi(z^0)$. Therefore,
$\liminf\Phi(z^s)\ge\liminf\Phi(\tilde z^s)=\Phi(z^0)$ that proves
lower semicontinuity of $\Phi$ at $z^0$ and thus its continuity.

{\it (b)} For any $t\in\Rn$ and $R>0$,
\begin{eqnarray*}
R^{-1}g'(\Phi,x,Rt) &= &R^{-1}\Phi(|x_1|+e^{Rt_1},\ldots,
|x_n|+e^{Rt_n})\\
&=& \Phi([|x_1|+e^{Rt_1}]^{1/R},\ldots,
[|x_n|+e^{Rt_n}]^{1/R}).
\end{eqnarray*}
The argument $[|x_k|+e^{Rt_k}]^{1/R}$ tends to $\exp\{t_k^+\}$ if
$x_k\neq 0$, and to $\exp\{t_k\}$ otherwise, so the statement follows
from {\it (a)}.

\medskip
The growth characteristics (\ref{eq:type}), (\ref{eq:mtype}),
(\ref{eq:dir}) of functions $u\in\LL$ can be expressed in terms of
the indicators. We will use the following notation:
\beq
\label{eq:ik}
\Ed=(1,\ldots,1),\ \Ed_1=(1,0,\ldots,0),\ \Ed_2=(0,1,0,\ldots,0),
\ \ldots, \Ed_n=(0,\ldots,0,1).
\eeq

\begin{prop}
\label{prop:char}
\begin{enumerate}
\item[(a)]
$\nui=\nu(\Psux,\infty)=\psux(\Ed)$;
\item[(b)]
$\nuai=\nu(\Psux,a,\infty)=\psux(a)\quad\forall a\in\Rnp$;
\item[(c)]
$\sigma_k(u)=\sigma_k(\Psux)=\psux(\Ed_k),\ k=1,\ldots,n$.
\end{enumerate}
\end{prop}

{\it Proof.}
The relation $\nuai=\psux(a)$ follows directly from the definition of
$\psux$. The equalities $\nui=\psux(\Ed)$ and
$\sigma_k(u)=\psux(\Ed_k)$ are proved in Theorems 1 and 2
of \cite{keyLe5}. Being applied to the function $\Psux$ instead of
$u$, they give us the first equalities in {\it (a)--(c)} in view of
Proposition \ref{prop:indind}.
The proof is complete.

\medskip
\beth
\label{theo:maj}
Let $u\in\LL$, $x\in\Cn$. Then
\beq
u(z)\le\Psux(z-x)+C\quad\forall z\in\Cn
\label{eq:maj}
\eeq
with $C=g'(u,x,0)$. Moreover, $\Psux$ is the least indicator
satisfying (\ref{eq:maj}) with some constant $C$.
\eth

{\it Proof.}
By Proposition \ref{prop:psi},
$$
g'(u,x,Rt)\le R\,\psux(t)+g'(u,x,0)\quad\forall R>0,\ \forall t\in\Rn,
$$
that implies (\ref{eq:maj}) since $u(z+x)\le g'(u,x,\cL(z))$.

If $\Phi\in\II$ satisfies $u(z)\le\Phi(z-x)+C$, then
$$
\Psux\le\Psi_{\Phi(\cdot+x),x}=\Psi_{\Phi,0}=\Phi,
$$
the latter equality being a consequence of Proposition
\ref{prop:indind}. The theorem is proved.

\medskip
The indicator $\Psux$ can be easily calculated in the algebraic case,
i.e. when $u$ is the logarihm of modulus of a polynomial mapping.
Recall that {\it the index} $I(P,x,a)$ of a polynomial $P$ at
$x\in\Cn$ with respect to the weight $a\in\Rnp$ is defined as
$$
I(P,x,a)=\inf\,\{\langle a,J\rangle :\: J\in\omega_x\}
$$
where
$$
\omega_x=\{J\in\Znp:\: {\partial^JP\over\partial z^J}(x)\neq 0\}
$$
(see e.g. \cite{keyLa}). For any $t\in\Rn$ we define
\beq
\label{eq:upind}
I_{up}(P,x,t)=\sup \,\{\langle t,J\rangle :\: J\in\omega_x\},
\eeq
the {\it upper index} of $P$ at $x\in\Cn$ with respect to $t\in\Rn$.
Clearly, $I_{up}(P,x,t)=-I(P,x,-t)$ for all $t\in -\Rnp$.

\begin{prop}
\label{prop:index}
Let $u=\log|P|,\ P:\Cn\to{\bf C}$ being a polynomial. Then
$$
\psux(t)= I_{up}(P,x,t) \quad\forall t\in\Rn,\ \forall x\in\Cn.
$$
\end{prop}

{\it Proof.}
Let
$$
P(z)=\sum_{J\in\omega_x}c_J\,(z-x)^J
$$
and $d=I_{up}(P,x,t)$, so $b_J:=\langle t,J\rangle -d\le 0$ $\forall
J\in\omega_x$. Then
$$
R^{-1}g'(u,x,t)=d+R^{-1}\sup_\theta\,\{\log|\sum_J
c_J\exp [Rb_Ji\langle\theta,J\rangle]|\}.
$$
Since there exists  $J_0\in\omega_x$ with
$b_{J_0}=0$, the second term here tends to $0$ as $R\to +\infty$, and
the statement follows.

\medskip
\begin{prop}
\label{prop:sup}
Let $u_1,\ldots,u_m\in\LL,\ u=\sup_k\,u_k,\ v=\log\sum_k\,\exp
u_k$. Then
$$
\Psux=\Psi_{v,x}=\sup_k\,\Psi_{u_k,x}.
$$
\end{prop}

{\it Proof.}
Since $u\ge u_k$, we have $\Psux\ge \sup_k\,\Psi_{u_k,x}$.
On the other hand, by (\ref{eq:maj}),
$$
u(z)\le \sup_k\,\{\Psi_{u_k,x}+C_k\}\le \sup_k\,\Psi_{u_k,x} +
\sup_k\,C_k,
$$
and the equality $\Psux=\sup_k\,\Psi_{u_k,x}$ results from
Theorem \ref{theo:maj}.

Similarly, the relations $\Psi_{v,x}\ge\Psux$ and
$$
v(z)\le \log\sum_k\,\exp [\Psi_{u_k,x}(z-x)+C_k]
\le\Psux(z-x)+m+ \sup_k\,C_k
$$
imply $\Psux=\Psi_{v,x}$, and the proof is complete.

\medskip
As a corollary of Propositions \ref{prop:index} and
\ref{prop:sup} we get

\begin{prop}
\label{prop:indmap}
Let
$$
u={1\over q}\log\sum_{k=1}^m |P_k|^q
$$
with $P_1,\ldots, P_m$ polynomials and $q>0$. Then
$\psux(t)=\sup_k\,I_{up}(P_k,x,t)$.
\end{prop}

\medskip
The indicator $\Psux$ can be described as a tangent (in
logarithmic coordinates) to the original function $u\in\LL$.
For $z\in\Cn$ and $m\in{\bf N}$, we set $z^m=(z_1^m,\ldots,z_n^m)$
and define the function
$$ (\tmux)(z)=m^{-1}u(x+z^m)\in\LL.
$$

\beth
\label{theo:tangent}
$\tmux \to \Psux$ in $L_{loc}^1(\Cn)$ as $m\to+\infty$.
\eth

{\it Proof.}
First, the family $\{\tmux\}_m$ is relatively compact in
$L_{loc}^1(\Cn)$. Really, (\ref{eq:maj}) implies
\beq
(\tmux)(z)\le\Psux(z-x)+m^{-1}C\quad\forall m.
\label{eq:t1}
\eeq
Therefore, the family is uniformly bounded above on each
compact subset of $\Cn$. Besides, $g(\tmux,0,0)=m^{-1}g(u,x,0)
\to 0$ and hence $g(\tmux,0,0)\ge -1$ for all $m\ge m_0$, and
the compactness follows.

Now let $v$ be a partial weak limit of $\tmux$, i.e.
$\tsmux\to v$ for some subsequence $m_s$. By (\ref{eq:t1}),
\beq
\label{eq:t2}
v\le\Psux.
\eeq
On the other hand, the convergence of $\tsmux$ to $v$ implies
$$
g(\tsmux,0,t)\to g(v,0,t)\quad\forall t\in\Rn.
$$
At the same time, by the definition of $\psux$,
$$
g(\tsmux,0,t)=m^{-1}g(u,x,mt)\to\psux(t)\quad\forall t\in\Rn,
$$
so $g(v,0,t)=\psux(t)$ and thus $g(v,0,t)=g(\Psux,0,t)$.
Being compared to (\ref{eq:t2}) it gives us $v=\Psux$, that
completes the proof.

\medskip
We conclude this section by studying dependence of $\Psux$ on
$x$.

\begin{prop}
\label{prop:generic}
Let $u\in\LL$. Then
\begin{enumerate}
\item[(a)]
$\Psi_u(z):=\sup\,\{\Psux(z):\:x\in\Cn\}\in\LL$;
\item[(b)]
for any $z\in\Cn$, $\Psux(z)=\Psi_u(z)$ for all
$x\in\Cn\setminus E_z$, $E_z$ being a pluripolar subset of
$\Cn$;
\item[(c)]
for any $z\in D^{-1}$, $\Psux(z)=\Psi_u(z)$ for all
$x\in\Cn$;
\item[(d)]
$\Psi_u(z)\ge 0\quad \forall z\in\Cn,\quad \Psi_u\equiv 0$ in
$D$.
\end{enumerate}
\end{prop}

{\it Proof.}
Since $u\in\LL$, there is a constant $A>0$ such that
$u(z)\le A\,S^+(z)\ \forall z\in\Cn$, where $S^+(z)
=S_1^+(z)=\sup_k\,\log^+|z_k|$.

We fix a point $z\in\Cn$ and consider the function
$$
u_R(x)=R^{-1}g'(u,x,R\,\cL(z)),\quad R>0.
$$
It is plurisubharmonic in $\Cn$, and
$$
u_R(x)\le R^{-1}A\,S^+(|x_1|+|z_1|^R,\ldots, |x_n|+|z_n|^R).
$$
Therefore, the family $\{u_R\}_{R>1}$ is uniformly bounded
above on compact subsets of $\Cn$, and
\beq
\label{eq:unif}
u_\infty(x):=\limsup_{R\to\ +\infty}u_R(x)\le
A\,S^+(z)\quad\forall x\in\Cn.
\eeq
Its regularization $u_\infty^*(x)=\limsup_{y\to x}u_\infty(y)$ is
then plurisubharmonic in $\Cn$ and bounded and so
$u_\infty^*=const$. We have $u_\infty (x)\le u_\infty^*(x)$ for
all $x\in\Cn$ with the equality outside a pluripolar set
$E_z\subset\Cn$.
We observe now that $u_\infty(x)=\Psux(z)$ and $u_\infty^*(x)=
\Psi_u(z)$, so {\it (b)} is proved.

Let $z^{(j)}\to z$, then the set
$$
E=\bigcup_{j=1}^\infty E_{z^{(j)}}\cup E_z
$$
is pluripolar. For $x\in\Cn\setminus E$,
$$
\Psi_u(z)=\Psux(z)=\lim_{j\to\infty}\Psux(z^{(j)})=
\lim_{j\to\infty}\Psi_u(z^{(j)}),
$$
that proves continuity of $\Psi_u$. Therefore,
$\Psi_u=\Psi_u^*$ is plurisubharmonic and belongs to $\LL$ in
view of (\ref{eq:unif}), that gives us {\it (a)}.

If $x,\ y\in\Cn$ and $a\in\Rnp$, we have for any $\epsilon>0$,
$g'(u,x,Ra)\le g'(u,y,(1+\epsilon)Ra)$ for all
$R>R_0(\epsilon,x,y)$, so $\psux(a)\le\psi_{u,y}(a)$ that
implies {\it (c)}.

Finally, {\it (d)} follows from the relation $\Psux|_D=0$
provided $\nu(u,x)=0$.

\section{Monge-Amp\`ere measures}

Now we pass to study the Monge-Amp\`ere measures of functions
$u\in\LL$. We can benefit by the plurisubharmonicity of the
growth characteristic $\Psi_u$ as well as by its specific
properties established in the previous section.

Any indicator $\Phi$ belongs to $L_{loc}^\infty(\Cstar)$, so
$(dd^c\Phi)^n$ is well defined on $\Cstar$. If
$\Phi\in L_{loc}^\infty(\Cn\setminus\{0\})$, then
$(dd^c\Phi)^n$ is defined on the whole space $\Cn$; the class
of such indicators will be denoted by $\II_0$.

Let $\TT$ denote the distinguished boundary $\{z\in\Cn:\:
|z_1|=\ldots=|z_n|=1\}$ of the unit polydisk $D$.

\begin{prop}
\label{prop:MAind}
Let $\Phi\in\II$. Then
\begin{enumerate}
\item[(a)]
$(dd^c\Phi)^n=0$ on $\Cstar\setminus \TT$;
\item[(b)]
if $\Phi\in\II_0$, then $(dd^c\Phi)^n=\tau_\Phi'\,\delta(0)+
\tau_\Phi''\, dm_T$ where $\tau_\Phi',\tau_\Phi''\ge 0$,
$\delta(0)$ is the Dirac measure at $0$, and $dm_T=(2\pi)^{-n}
d\theta_1\ldots d\theta_n$ is the normalized Lebesgue measure on
$\TT$.
\end{enumerate}
\end{prop}

{\it Proof.}
{\it (a)} It suffices to show that for every $y\in\Cstar$ there
exists an analytic disk $\gamma_y$ containing $y$ such that the
restriction of $\Phi$ to $\gamma_y$ is harmonic near $y$
(\cite{keyFS}, Lemma~6.9). Let $y=(|y_1|e^{i\theta_1},\ldots,
|y_n|e^{i\theta_n})\in\Cstar$. Consider the mapping $\lambda_y:
{\bf C}\to\Cn$ given by
$$
\lambda_y(\zeta)=(|y_1|^\zeta
e^{i\theta_1},\ldots,|y_n|^\zeta e^{i\theta_n});
$$
note that $\lambda_y(1)=y$.
Since $y\in\Cstar\setminus \TT$, $\lambda_y$ is not constant.
Set $\Delta=\{\zeta\in{\bf C}: |\zeta -1|<1/2\}$ and $\gamma_y=
\lambda_y(\Delta)\subset\Cstar$. Then
$\Phi(\lambda_y(\zeta))={\rm Re}\:\zeta\cdot\Phi(\lambda_y(1))$,
so the restriction of $\Phi$ to $\gamma_y$ is harmonic.

{\it (b)} follows from {\it (a)} since locally
plurisubharmonic functions cannot charge pluripolar sets and
$\Phi(y)$ is independent of ${\rm arg}\,y_k$, $\okn$.

\medskip
We will say that the unbounded locus of $u\in PSH(\Cn)$ is {\it
separated at infinity} if there exists an exhaustion of $\Cn$
by bounded pseudoconvex domains $\Omega_k$ such that $\inf\,
\{u(z):\, z\in\partial\Omega_k\}>-\infty$ for each $k$. The
collection of all functions $u\in\LL$ whose unbounded loci are
separated at infinity will be denoted by $\LL_*$. By
\cite{keyD1}, Corollary 2.3, the Monge-Amp\`ere current
$(dd^cu)^n$ is well defined on $\Cn$ for any
function $u\in\LL_*$. Note also that $\Psux\in\II_0\ \forall
x\in\Cn$ for any $u\in\LL_*$.

We are going to compare the total Monge-Amp\`ere mass
$$
\mucn=\int_{\Cn} (dd^cu)^n
$$
of $u\in\LL_*$ with that of its indicator. The key result is
the following comparison theorem (which is actually a variant of
B.A.~Taylor's theorem \cite{keyT}).

\beth
\label{theo:comp1}
Let $u,v\in\LL_*,\ v\ge 0$ outside a bounded set, and
$$
\limsup_{z\to\infty} {u(z)\over v(z)+\eta \log|z|}\le 1\quad
\forall \eta>0.
$$
Then $\mucn\le M(v;\Cn)$.
\end{theo}

{\it Proof.} By the definition of the class $\LL_*$, there
exist numbers $0\le m_1\le m_2\le\ldots$ such that $u(z)>-m_k$
near $\partial\Omega_k$, $\{\Omega_k\}$ being the pseudoconvex
exhaustion of $\Cn$. Let $w(z)=\sup\,\{u(z),-m_k\}$ for
$\Omega_k\setminus\Omega_{k-1}$ (assuming
$\Omega_0=\emptyset$), then $w\in\LL\cap L_{loc}^\infty(\Cn)$
and satisfies the same asymptotic relation at infinity as $u$
does.
Besides,
$$
\int_{\Omega_k} (dd^cu)^n=\int_{\Omega_k} (dd^cw)^n\quad\forall k,
$$
so
\beq
\label{eq:c1}
\mucn=M(w;\Cn).
\eeq
Denote $v_\eta(z)=v^+(z)+\eta\log^+|z|$,
$\eta>0$. Let $\epsilon>0$ and
$C>0$, then $w(z)\le (1+\epsilon)v_\eta-2C$ for all $z\in
\Cn\setminus B_\alpha$ with $B_\alpha$ a ball of the radius
$\alpha=\alpha(\eta,C,\epsilon)$. Therefore,
$$
E(\eta,C,\epsilon):=\{z\in\Cn:\: (1+\epsilon)v_\eta-C<w(z)\}
\Subset B_\alpha.
$$
By the Comparison Theorem for bounded plurisubharmonic
functions,
\beq
\label{eq:c2}
\int_{E(\eta,C,\epsilon)} (dd^cw)^n \le
\int_{E(\eta,C,\epsilon)} (dd^c[(1+\epsilon)v_\eta-C])^n
\le (1+\epsilon)\int_{\Cn} (dd^c v_\eta)^n.
\eeq
For any compact $K\subset\Cn$ one can find $C>0$ such that
$K\subset E(\eta,C,\epsilon)$, so (\ref{eq:c2}) gives us
$$
M(w;\Cn)\le (1+\epsilon)M(v_\eta;\Cn).
$$
Since $v_\eta$ decreases to $v^+$ as $\eta\to 0$ and in view of
the arbitrary choice of $\epsilon$, we then get
$$
M(w;\Cn)\le M(v^+;\Cn)=M(v;\Cn),
$$
which by (\ref{eq:c1}) completes the proof.

\medskip
As an immediate consequence we have
\beth
\label{theo:comp2}
For any $u\in\LL_*$,
$$
\mucn\le M(\Psux;\Cn)\le M(\Psi_u;\Cn).
$$
\eth

\medskip
To get effective bounds for $\mucn$, we estimate the Monge-Amp\`ere
masses of the indicators.

\begin{prop}
\label{prop:max}
Let $\Phi\in\II,\ z^0\in D^{-1}$. Then
$$
\Phi(z)\le\Phi(z^0)\,\sup_k\,{\log^+|z_k|\over\log|z_k^0|}\quad
\forall z\in\Cn.
$$
\end{prop}

{\it Proof.}
Denote $\psi=[\Phi(z^0)]^{-1}\cE^*\Phi,\ a=\cL(z^0)\in\Rnp$.
It suffices to prove the relation
$\psi(t)\le s_a^+(t)$ for all $t$ with $|t_k|<a_k$, $\okn$; here
$s_a=\cE^*S_a$ with $S_a$ defined in (\ref{eq:Sa}), so
$s_a^+(t)=\sup_k\,t_k^+/a_k$.

We fix such a point $t$ and denote $\alpha=s_a^+(t)<1$,
$\beta=(1-\alpha)^{-1}>1$. Consider the segment
$l_t=\{a+\lambda\,(t-a):\: 0\le \lambda\le\beta\}\subset\Rn$; observe that
\beq
\label{eq:segment}
a_k+\beta\,(t_k-a_k)=\beta\,a_k\,(t_k a_k^{-1}-\alpha)\le 0,\
\okn.
\eeq
Let $u(\lambda)$ be the restriction of $\psi$ to $l_t$,
$v(\lambda)=(\alpha-1)\lambda+1$. The function $u$ is convex on $l_t$,
and $v$ is linear. Besides, $u(0)=v(0)=1$, $v(\beta)=0$, and
$u(\beta)\le 0$ in view of (\ref{eq:segment}).
Therefore, $u\le v$ on $l_t$. In particular, $\psi(t)=u(1)\le v(1)
=s_a^+(t)$, and the proposition is proved.

\medskip
\beth
\label{theo:ocenka}
Let $\Phi\in\II_0,\ z^0\in D^{-1}$. Then
$$
M(\Phi;\Cn)\le {[\Phi(z^0)]^n\over \log|z_1^0|\ldots\log|z_n^0|}.
$$
In particular, for any $u\in\LL$,
$$
M(\Psi_u;\Cn)\le {[\nuai]^n\over a_1\ldots a_n}\quad\forall a\in\Rnp.
$$
\eth

{\it Proof.}
The first relation follows from Proposition \ref{prop:max} and
Theorem \ref{theo:comp1}, since (taking $a=\cL(z^0)\in\Rnp$)
$$
M(S_a^+;\Cn)=M(S_a;\Cn)=[\log|z_1^0|\ldots\log|z_n^0|]^{-1}.
$$

The second inequality results now from Proposition~\ref{prop:char} {\it (b)}.

\medskip
We can give a geometric interpretation for the total Monge-Amp\`ere
masses of indicators, which in many cases leads to their exact
calculation.

Let $\Phi\in\II,\ \varphi=\cE^*\Phi$. Denote
$$
\Theta_\Phi^+=\{a\in\Rn:\: \langle a,t\rangle\le\varphi^+(t)\quad
\forall t\in\Rn\}.
$$

\begin{prop}
\label{prop:gamma}
$\Theta_\Phi^+$ is a convex compact subset of $\overline{\Rnp}$,
$\Theta_\Phi^+\subset\{a\in\Rn:\: 0\le a_k\le\varphi^+(\Ed_k),
\ \okn\}$.
\end{prop}

{\it Proof.} Convexity of $\Theta_\Phi^+$ is evident. Further, if
$a\in\Theta_\Phi^+$, then $\langle a,\pm\Ed_k\rangle
\le\varphi^+(\pm\Ed_k)$, $\Ed_k$ being defined by (\ref{eq:ik}),
and the statement follows
because $\varphi^+(-\Ed_k)=0$.

\medskip
By $Vol(P)$ we denote the Eucledean volume of $P\subset\Rn$.

\beth
\label{theo:volume}
For any $\Phi\in\II_0$,
$$
M(\Phi;\Cn)=n!\,Vol(\Theta_\Phi^+).
$$
\eth

{\it Proof.}
By Proposition \ref{prop:MAind},
$$
M(\Phi;\Cn)=M(\Phi^+;\Cn)=\int_\TT (dd^c\Phi^+)^n.
$$

It can be easily checked that the complex Monge-Amp\`ere operator $(dd^cU)^n$
of an $n$-circled locally bounded plurisubharmonic function $U$
is related to the real Monge-Amp\`ere operator ${\cal MA}[u]$
of its convex image $u$ by the equation
$$
\int_G(dd^cU)^n=n!\,\int_{\cL(G)}{\cal MA}[u]
$$
for every $n$-circled Borelean set $G\subset\Cn$ (see e.g.
\cite{keyR}).
Since $\cL(\TT)=\{0\}$,
\beq
M(\Phi^+;\Cn)=n!\,{\cal MA}[\varphi^+](\{0\}).
\label{eq:v1}
\eeq
As was established in \cite{keyRaT}, for any convex function $v$
on a domain $\Omega\subset\Rn$,
\beq
{\cal MA}[v](F)=Vol\,(\omega(F,v))\quad\forall F\subset\Omega,
\label{eq:v2}
\eeq
where
$$
\omega(F,v)=\bigcup_{t^0\in F}\{a\in\Rn:\: v(t)\ge v(t^0)+
\langle a, t-t^0\rangle\ \forall t\in\Omega\}
$$
is the gradient image of the set $F$ for the surface
$\{y=v(x),\ x\in\Omega\}$. In our situation, it means that
\beq
\omega(\{0\},\varphi^+)=\{a\in\Rn:\: \varphi^+(t)\ge
\langle a,t\rangle\ \forall t\in\Rn\}=\Theta_\Phi^+,
\label{eq:v3}
\eeq
so the statement follows from (\ref{eq:v1})--(\ref{eq:v3}).

\medskip
The set $\Theta_\Phi^+$ for $\Phi=\Psux,\ u\in\LL,\ x\in\Cn$, will be
denoted by $\Theta_{u,x}^+$, and by $\Theta_u^+$ for $\Phi=\Psi_u$.
Then Theorems \ref{theo:comp2} and \ref{theo:volume} give us

\beth
\label{theo:Kou}
For any $u\in\LL_*$ and $x\in\Cn$,
\beq
\label{eq:Kou}
\mucn\le n!\,Vol(\Theta_{u,x}^+)\le n!\,Vol(\Theta_{u}^+).
\eeq
\eth

{\it Remark.} Let $u=\log|P|,\ P=(P_1,\ldots, P_N)$ being a
polynomial mapping. By Proposition \ref{prop:indmap},
$$
\psux(t)=I_{up}(P,x,t)=\sup_{1\le j\le N}\, I_{up}(P_j,x,t),
$$
the upper indices $I_{up}(P_j,x,t)$ defined by (\ref{eq:upind}).
In this case,
$$
\Theta_{u,x}^+=\{a\in\Rn:\: \langle a,t\rangle\le
I_{up}^+(P,x,t)\ \forall t\in\Rn\},
$$
so $\Theta_{u,0}^+$ coincides with the Newton polyhedron for $P$ at
infinity (see Introduction). If $N=n$ and $P^{-1}(0)$ is discrete,
then $\mucn$ is the number of zeros of $P$ counted with the
multiplicities. For this case, (\ref{eq:Kou}) gives the bound
due to Kouchnirenko \cite{keyKou}.

\medskip
Theorem \ref{theo:Kou} produces also an upper bound for $\mucn$ via
the multitype $(\sigma_1(u),\ldots,\sigma_n(u))$ of the function $u$:

\beth
\label{theo:Dyson}
Let $u\in\LL_*$, then
\beq
\label{eq:Dys1}
\mucn\le n!\, \sigma_1(u)\ldots\sigma_n(u);
\eeq
in particular,
\beq
\label{eq:Dys2}
\sum_x[\nux]^n\le n!\,\sigma_1(u)\ldots\sigma_n(u).
\eeq
\eth

{\it Proof.}
By Propositions \ref{prop:char} and \ref{prop:gamma},
$\Theta_{u}^+\subset\{a\in\Rn:\: 0\le a_k\le\sigma_k(u),\ \okn\}$,
and (\ref{eq:Dys1}) follows from Theorem \ref{theo:Kou}. It implies
(\ref{eq:Dys2}) in view of the known inequality
$(dd^cu)^n|_{\{x\}}\ge [\nux]^n$.

\medskip
{\it Remark.} It can be shown that inequality (\ref{eq:Dys1})
implies Dyson's lemma for algebraic hypersurfaces with isolated singular
points (see \cite{keyVi}).

\vskip 0.5cm

Mathematical Division, Institute for Low Temperature Physics
\par
47 Lenin Ave., Kharkov 310164,
Ukraine

\vskip0.1cm

E-mail: \quad rashkovskii@ilt.kharkov.ua, 
\quad rashkovs@rashko.ilt.kharkov.ua



\end{document}